\documentclass[12pt]{article}
\usepackage{amsmath}
\usepackage{amsfonts}
\title{Genus Zero Modular Functions}
\author{Bong H. Lian and Joshua L. Wiczer}
\date{31 Aug 2003}
\newtheorem{theorem}{Theorem}
\begin{document}
\maketitle

This project was sponsored through the Schiff Fellowship program of
Brandeis University. This project involved using the power series
method to construct a third order nonlinear ordinary differential equation, 
a Schwarzian equation, for each of
the ``genus zero" modular functions, described in the Conway-Norton paper.  
We first use the Borcherd recursion formuli to
generate, in each case, a modular function up to whatever degree we desire, and
then use the fact that there is a Schwarzian equation, 
determined by a single rational function we call a Q-value. By similar
power series methods, we compute the coefficients of our rational
function, and hence have all the necessary data to create a Schwarzian
differential equation for each modular function. This equation can, in turn, be used to
recover the modular function itself.

\section{Basic Notions of Differential Equations}
First, we look at some basic notions of differential equations,
that will form the basis for later work.  One of the most powerful
theorems is the existence and uniqueness theorem for differential
equations.

\begin{theorem}
If $\dfrac{b(t)}{a(t)}$ and $\dfrac{c(t)}{a(t)}$ are continuous
functions on an open interval $I$ that includes $t=t_{0}$, then
the initial value problem
$$a(t)y''(t)+b(t)y'(t)+c(t)y(t)=0,$$ $$y(t_{0})=a,$$ $$y'(t_{0})=b$$
has a unique solution on $I$.
\end{theorem}

So, if we get a differential equation of that form, we will be
able to find a unique solution.  One way to find a solution is to
use the ``method of Frobenius."  We assume the solution takes on
the form $$y(t)=x^{r}\sum_{n=0}^{\infty}a_{n}t^{n}=
\sum_{n=0}^{\infty}a_{n}t^{r+n}.$$  If $y(t)$ is of that form, then
$y'(t)$ is of the form
$$y'(t)=\sum_{n=0}^{\infty}a_{n}(r+n)t^{r+n-1},$$
and then $y''(t)$ is of the from
$$y''(t)=\sum_{n=0}^{\infty}a_{n}(r+n)(r+n-1)t^{r+n-2}.$$
Then we plug those values into our second-order differential
equation, and find as many coefficients as we desire (by solving a
system of equations).

\section{Hypergeometric Functions}
If we have an ODE of the form
$$\dfrac{d^{2}u}{dz^{2}}+P(z)\dfrac{du}{dz}+Q(z)u=0$$ and $P$ and
$Q$ are rational functions of the form
$\dfrac{\sum_{n=0}^{i}\alpha_{n}z^{n}}{\sum_{n=0}^{j}\beta_{n}z^{n}}$,
with $P$ and $Q$ having no common factors, and if $P$ and $Q$ are
regular except possibly at the poles 0, 1, or $\infty$, then we can
do a linear fraction transformation to 0, 1, or $\infty$, giving us
an ODE of the form
$$z(1-z)\dfrac{d^{2}u}{dz^{2}}+[c-(a+b+1)z]\dfrac{du}{dz}-abu=0,$$
where a, b, c $\in \mathbb{C}$.  Details on the precise nature of
the transformation can be found in Seaborn.  This equation is called the \emph{Gauss hypergeometric equation}.

Now, to find a solution to this ODE, we use the method of
Frobenius, and say the solution is of the form
$$u(z)=\sum_{n=0}^{\infty}a_{n}z^{n+s}$$ and can solve to get a recursive definition
of the coefficients. More specifically, we get the recursion
$$a_{n+1}=\frac{(n+a)(n+b)}{(n+1)(n+c)}a_{n}.$$

Following the recursion through, a bit, we see that
$$\begin{array}{lll}a_{1} & = & \dfrac{ab}{c}a_{0} \\
a_{2} & = & \dfrac{a(a+1)b(b+1)}{2c(c+1)}a_{0} \\
a_{3} & = & \dfrac{a(a+1)(a+2)b(b+1)(b+2)}{1\cdot 2\cdot 3
c(c+1)(c+2)}a_{0} \\
\vdots & & \vdots \\
a_{n} & = & \dfrac{a(a+1)(a+2)\cdots (a+n-1)b(b+1)(b+2)\cdots
(b+n-1)}{n! c (c+1)(c+2) \cdots (c+n-1)}a_{0} \\ \end{array}$$

As a notational convention, we define the \emph{Pochhamer symbol},
$(a)_{n}$ as:
$$\begin{array}{lll}(a)_{0}& =&1 \\ (a)_{n} & = & a(a+1)(a+2)\cdots
(a+n-1) \end{array}$$ and we can write the solution to a
differential equation of that form as
$$u(z)=\sum_{n=0}^{\infty}
\dfrac{(a)_{n}(b)_{n}}{n!(c)_{n}}a_{0}z^{n},$$ which we call a \emph{hypergeometric function}.

\section{The Schwarzian}
The Schwarzian derivative is a ``differential operator invariant
under the group of all fractional linear transformation" (Hille).
The Schwarzian derivative is defined by
$$\{w,z\}=\frac{w'''}{w'}-\frac{3}{2}\left(\frac{w''}{w'}
\right)^{2},$$ where $w$ is dependent on $z$.

Now we have our key theorem, taken from Hille:
\begin{theorem}  Let $y_{1}$ and $y_{2}$ be two linearly
independent solutions of the equation $y''+Q(z)y=0$, and be
defined and analytic is a simply connected domain D in the complex
plane.  Then $w(z)=\frac{y_{1}(z)}{y_{2}(z)}$ satisfies
$\{w,z\}=2Q(z)$ at all points of D where $y_{2}(z)\neq 0$.
Conversely, if $w(z)$ is a solution of $\{w,z\}=2Q(z)$, analytic
in some neighborhood of a point $z_{0}\in D$, then you can find
linearly independent solutions, $u(z)$ and $v(z)$ (of
$y''+Q(z)y=0)$ defined in D such that $w(z)=\frac{u(z)}{v(z)}$,
and if $v(z_{0})=1$, then solutions $u$ and $v$ are uniquely
defined.
\end{theorem}

We can find $\{z,w\}=2Q(z)x'(z)^{2}$, the interchanging of dependent and independent variable, by primarily using
the inverse function theorem, stating that if $g$ is the inverse of $f$, then $g'(z)=\frac{1}{f'(g(z))}$.

We can also transform a Gauss hypergeometric equation into an equation of the form $y''(z)+Q(z)y(z)=0$, and therefore
we can see that every Gauss hypergeometric equation will have a corresponding ``Q"-value.  In particular,
we can compute that the ``Q"-value will be
$$Q=\dfrac{-c^{2}+z(2ab(z-2)+z-a^{2}z-b^{2}z)+2c(1+(a+b-1)z)}{4(z-1)^{2}z^{2}}.$$
There will be singularities at $z=0$, $z=1$, and at $\infty$.

\section{The Picard-Fuchs Equation}
As an example, let us consider the differential equation
$$\dfrac{d^{2}\Omega}{dJ^{2}}+\dfrac{1}{J}\frac{d\Omega}{dJ}+\dfrac{31J-4}{144J^{2}(1-J)^{2}}\Omega=0.$$
It has two linearly independent solutions:
$$\begin{array}{lll} f_{1} & = & 1 + a_{1}J+a_{2}J^{2}+\ldots =
1+O(J) \\
f_{2} & = & f_{1}\log{J}+b_{1}J+b_{2}J^{2}+\ldots = f_{1}+ \log{J} + O(J)
\\ \end{array}$$

We can normalize the ODE, by letting $\widetilde{f}=gf$, giving us

$$\frac{d^{2}\widetilde{f}}{dJ^{2}}+Q \widetilde{f}=0.$$

We are going to examine $t=\dfrac{f_{2}}{f_{1}}$, and then, we will
look at $q=e^{t}$, which will be of the form
$$ J + c_{2}J^{2}+c_{3}J^{3}+\ldots .$$
This is because
$$\dfrac{f_{2}}{f_{1}}=\dfrac{f_{1}\log{J}+O(J)}{f_{1}}=\log{J}+O(J),$$
making $q=e^{\frac{f_{2}}{f_{1}}}=J+e^{O(J)} = O(J)$.  Moreso,
hence
$$q=z+a_{2}z^{2}+a_{3}z^{3}+\ldots.$$

Going back to the original Picard-Fuchs equation, we clear all
denominators, and normalize it to turn it into
$$\widetilde{f}''(q)+\frac{1-1968q+2654208q^{2}}{4(1-1728q)^{2}q^{2}}\widetilde{f}(q).$$
We can read off its $Q$ value as
$$\frac{1-1968q+2654208q^{2}}{4(1-1728q)^{2}q^{2}}.$$

Looking at this from another angle, we can find the solution to
the Picard-Fuchs equation by ansatz.  To find the first solution,
we let
$$\begin{array}{lll}
f_{1}(q) &=&\sum_{n=0}^{\infty}a_{n}q^{r+n} \\
f_{1}'(q) &=&\sum_{n=0}^{\infty}a_{n}(r+n)q^{r+n-1} \\
f_{1}''(q)&=&\sum_{n=0}^{\infty}a_{n}(r+n)(r+n-1)q^{r+n-2}
\end{array}$$
Plugging everything in, and letting $a_{0}=1$, we can recursively
find the rest of the coefficients, giving us
$$f_{1}(q)=1-372q-266076q^{2}-277702608q^{3}-336151952604q^{4}+O(q)^{5}.$$

Our second solution will be of the form
$$f_{2}(q)=f_{1}\log{q}+\sum_{n=0}^{\infty}b_{n}q^{n}.$$
By similar plugging in, and letting $b_{0}=0$, we get our second
solution to look like
$$f_{2}(q)=f_{1}\log{q}+744q+196884q^{2}+77574992q^{3}+9847760106q^{4}+O(q)^{5}.$$

We can now find the inverse of this, which gives us
$$ q - 744q^{2}+356652q^{3}-140361152q^{4} +O(q)^{5} $$
and taking one-over this, we get
$$ \frac{1}{q}+744+196884q-104836240q^{2}+O(q)^{3}.$$

This is called the elliptic modular function $j$.

\section{Modular Functions and Their Q-values}

Let $h$ be the complex upper-half plane.  A modular function is a holomorphic function $f:h \rightarrow
\mathbb{C}$ if it is invariant under the group $SL(2,\mathbb{Z})$.  $SL(2,\mathbb{Z})$ is
defined as all matrices of the form

$$ \left\{ \left[\begin{array}{ll}a&b \\ c&d \end{array}\right] \left|
\begin{array}{l} a,b,c,d \in \mathbb{Z} \\ ad-bc=1 \end{array} \right. \right\}.$$
$f$ is said to be invariant under a matrix $A$ if $f(A(t)) = f(t)$, where $A(t) = \frac{at+b}{ct+d}$.

 There are 175 so-called ``genus zero" modular functions, of which $j(q)$ is one.  The genus-zero modular functions arise in the work of Conway and Norton, in their study of the Moonshine Conjecture.  It is known that each of these functions has a $q$-series expansion of the form $q^{-1}+c+O[q]$.

Conway and Norton give us the first 12
coefficients of these modular functions.  For example, Conway and Norton give us the values:
1, 24, 196884, 21493760, 864299970, 20245856256, 333202640600, 4252023300096, 44656994071935, 401490886656000, 3176440229784420, and 22567393309693600.

These values correspond to the modular function
\begin{eqnarray*}
&&q^{-1} + 24 + 196884q + 21493760q^{2} + 864299970q^{3} + 20245856256^{4} + 333202640600q^{5}\\&& + 4252023300096q^{6}
+ 44656994071935q^{7} + 401490886656000q^{8} \\&&+ 3176440229784420q^{9} + 22567393309693600q^{10}+\cdots
\end{eqnarray*}


However, it is important to understand that the general form of a modular function transformed into a modular function of form 1A is

$$q^{-1} + c + 196884q + 21493760q^2 + ...$$
where $c\in \mathbb{R}$.  Because of this, we modified Conway and
Norton's data, to let all of the constants (coefficients of $q^{0}$)
become 0.  This will make our Q-values look nicer, and there's a
transform we can use to find the Q-value for any of the 175
geneus-zero modular functions, with any constant coefficient.

Borcherds provides us with a recursive algorithm to determine as many coefficients as we
want. So now, our task is to determine the unique ``Q" value for
each of the 175 genus-zero modular functions, with the constant-coefficient set at 0.

All of Conway and Norton's coefficients have been inputted into a text file, \texttt{Coeff.txt}.  (All of these files are available online at people.brandeis.edu/lian/)  The format of the data is as follows:

\begin{itemize}
\item Each set of coefficients, given by Conway and Norton as 1A, 2A, 2B, 3A, etc. has been numbered from 1 to 175.  The function \texttt{h} inputs the Conway and Norton numeration and outputs my numeration.  For example, \texttt{h[1A]=1} and \texttt{h[10E]=33}.

\item The coefficients of $q^{-1}$ up through $q^{10}$, as listed in the Conway and Norton paper, are provided in the two variable function \texttt{d}.  The format for \texttt{d} is \texttt{d[}number of the modular function, degree of $q$] = coefficient of $q$.  So, for example, \texttt{d[1,-1]=1}, and \texttt{d[24,1]=-4}.

\item Borcherds' recursion formuli for $c_{g}(k)$ depend on values of $c_{g}(j)$ where $j<k$, but they also depend on values of $c_{g^{2}}(j)$, where $j<k$.  For each set of coefficients, $c_{g}$, there is a corresponding set of coefficients, $c_{g^{2}}$.  These are given (albeit a bit cryptically) in Conway and Norton's paper (in Table 2), but can also be determined by a trial-and-error experiment, which performed.  (Breifly, the experiment involved starting with $c_{g}$ being the coefficients for 1A, and then running through all the other 174 modular functions, to see if when those coefficients were plugged in for $c_{g^{2}}$, we would get accurate coefficients for $c_{g}$.)  To summarize the findings briefly, for the cofficients for the modular function $XY$, where $X$ is a number between 1 and 119 (once again, corresponding to the numberings of the modular functions by Conway and Norton), and where Y is a letter.  If $X$ is even, then the corresponding $c_{g^{2}}$ is $(X/2)Z$, where $Z$ is some letter (not predictable).  If $X$ is odd, then the corresponding $c_{g^{2}}$ is $XY$.  The data is stored in the function $f$.  The format for $f$ is $f$[number of the modular function]=number of the modular function for $c_{g^{2}}$.
\end{itemize}

We wrote a Mathematica program (\texttt{ModularFn.txt}) that compues any Conway and Norton's delineated 175 genus-zero modular functions up to any given degree.  We first input a value \texttt{md} which correspondings to the Conway and Norton name of the modular function.  We then input a value \texttt{n} which corresponds to the degree of the modular function we wish to compute.  The program then uses the Bourcherds recursion formuli to determine all the coefficients up to \texttt{n}.

Once we have a modular function, we can compute the Q-value.  A Q-value determines a modarul function, and hence encodes the same information as the function itself.  We know that Q will be a rational function, with finitely many terms in the numerator and denominator.  To compute the Q-values, we wrote a program in Mathematica, \texttt{Qvalues.txt}, that inputs a modular function, \texttt{md}, inputs a maximum degree of the numerator, \texttt{r}, and inputs a maximum degree of the denomninator \texttt{s}.  After we compute the modular function to a degree of 300, using the same basic structure as the program \texttt{ModularFn.txt}.  We then return to the Schwarzian equation, plugging in our modular functions and the Q-value

$$Q=\frac{\sum_{i=0}^{r}b_{i}z(q)^{i}}{\sum_{i=0}^{s}c_{i}z(q)^{i}}.$$
We can bring everything to one side of the equation, and solve for
set all the coefficients of all the ``q"-values equal to 0:
$q^{2}, q^{3}, q^{4}$ and so on, up to $s+r+2$. If Mathematica isn't
able to solve the equations (to get anything other than 0 for all
the $b_{i}$s and $c_{i}$s, then we need to increase the degree of
the numerator and denominator, and keep increasing them until
we are able to get good values for them all.

We carried this out for all 175-modular functions to give us the
following Q-values:

\begin{tabbing}
Number \hspace{.25in} \= Q-value \hspace{.1in} \= \\

1A \> $\frac{1 - 480z + 1743552z^{2}}{4z^{2}(-1 + 240z + 732096z^{2})^{2}}$ \\
\\
2A \> $\frac{1 - 96z + 40640z^{2}}{4z^{2}(-1 + 48z + 15808z^{2})^{2}}$ \\
\\
2B \> $\frac{1 + 32z + 2752z^{2}}{4z^{2}(1 + 16z - 960z^{2})^{2}}$ \\
\\
3A \> $\frac{1 - 48z + 7560z^{2}}{4z^{2}(-1 + 24z + 2772z^{2})^{2}}$ \\
\\
3B \> $\frac{1 + 6z + 513z^{2}}{4z^{2}(1 + 3z - 180z^{2})^{2}}$ \\
\\
3C \> $\frac{1 + 9936z^{3}}{4z^{2}(1 - 1728z^{3})^{2}}$ \\
\\
4A \> $\frac{1 - 32z + 2752z^{2}}{4z^{2}(-1 + 16z + 960z^{2})^{2}}$ \\
\\
4B \> $\frac{1 + 320z^{2} + 49152z^{4}}{4z^{2}(1 - 256z^{2})^{2}}$\\
\\
4C \> $\frac{1 + 192z^{2}}{4z^{2}(1 - 64z^{2})^{2}}$\\
\\
4D \> $\frac{1 - 64z^{2} + 4096z^{4}}{4z^{2}(1 + 64z^{2})^{2}}$\\
\\
5A \> $\frac{1 + 8z + 848z^{2} + 38096z^{3} + 359872z^{4}}{4z^{2}(1 + 4z - 656z^{2} - 7424z^{3})^{2}}$\\
\\
5B \> $\frac{1 + 8z + 98z^{2} + 1096z^{3} + 2497z^{4}}{4z^{2}(1 + 4z - 31z^{2} - 174z^{3})^{2}}$\\
\\
6A \> $\frac{1 + 4z + 492z^{2} + 16352z^{3} + 120608z^{4}}{4z^{2}(1 + 2z - 388z^{2} - 3080z^{3})^{2}}$\\
\\
6B \> $\frac{1 + 4z + 486z^{2} + 16804z^{3} + 117841z^{4}}{4z^2(1 + 2z - 383z^{2} - 3156z^{3})^{2}}$\\
\\
6C \> $\frac{1 + 4z + 108z^{2} - 800z^{3} + 2848z^{4}}{4z^{2}(1 + 2z - 68z^{2} + 120z^{3})^{2}}$\\
\\
6D \> $\frac{1 + 4z + 6z^{2} + 1124z^{3} + 1457z^{4}}{4z^{2}(1 + 2z + 17z^{2} - 164z^{3})^{2}}$\\
\\
6E \> $\frac{1 + 4z + 54z^{2} + 388z^{3} + 769z^{4}}{4z^{2}(1 + 2z - 23z^{2} - 60z^{3})^{2}}$\\
\\
6F \> $\frac{1 - 304z^{3} + 4096z^{6}}{4z^{2}(1 + 64z^{4})^{2}}$\\
\\
7A \> $\frac{1 + 2z + 313z^{2} + 8904z^{3} + 50328z^{4}}{4z^{2}(1 + z - 252z^{2} - 1620z^{3})^{2}}$\\
\\
7B \> $\frac{1 + 2z + 19z^{2} + 378z^{3} + 201z^{4}}{4z^{2}(1 + z - 7z^{2} - 52z^{3})^{2}}$\\
\\
8A \> $\frac{1 + 224z^{2} + 5120z^{3} + 22784z^{4}}{4z^{2}(-1 + 176z^{2} + 896z^{3})^{2}}$\\
\\
8B \> $\frac{1 + 64z^{2} + 4096z^{4}}{4z^{2}(1 - 64z^{2})^{2}}$\\
\\
8C \> $\frac{1 + 2816z^{4}}{4z^{2}(1 - 256z^{4})^{2}}$\\
\\
8D \> $\frac{(1 - 16z^{2})^{2}}{4z^{2}(1 + 16z^2)^{2}}$\\
\\
8E \> $\frac{(1 + 16z^{2})^{2}}{4z^{2}(1 - 16z^{2})^{2}}$\\
\\
8F \> $\frac{1 - 640z^{4} + 4096z^{8}}{4z^{2}(1 + 64z^{4})^{2}}$\\
\\
9A \> $\frac{1 + 162z^{2} + 3456z^{3} + 14337z^{4}}{4z^{2}(-1 + 135z^{2} + 594z^{3})^{2}}$\\
\\
9B \> $\frac{1 + 216z^{3}}{4z^{2}(1 - 27z^{3})^{2}}$\\
\\
10A \> $\frac{1 + 128z^{2} + 2256z^{3} + 10560z^{4}}{4z^{2}(-1 + 112z^{2} + 384z^{3})^{2}}$\\
\\
10B \> $\frac{1 + 32z^{2} - 304z^{3} + 1088z^{4} - 2304z^{5} + 3072z^{6}}{4z^{2}(1 - 32z^{2} + 64z^{3})^{2}}$\\
\\
10C \> $\frac{1 - 22z^{2} + 256z^{3} + 185z^{4}}{4z^{2}(1 + 13z^{2} - 34z^{3})^{2}}$\\
\\
10D \> $\frac{1 + 10z + 147z^{2} + 3716z^{3} + 37683z^{4} + 159706z^{5} + 244897z^{6}}{4z^{2}(-1 - 5z + 107z^{2} + 961z^{3} + 2130z^{4})^{2}}$\\
\\
10E \> $\frac{1 + 4z + 16z^{2} + 124z^{3} + 630z^{4} + 1292z^{5} + 1408z^{6} + 1172z^{7} + 793z^{8}}{4z^{2}(-1 - 2z + 2z^{2} + 20z^{3} + 47z^{4} + 30z^{5})^{2}}$\\
\\
11A \> $\frac{1 + 8z + 112z^{2} + 2620z^{3} + 23424z^{4} + 86896z^{5} + 115884z^{6}}{4z^{2}(-1 - 4z + 88z^{2} + 668z^{3} + 1272z^{4})^{2}}$\\
\\
12A \> $\frac{1 - 4z + 108z^{2} + 800z^{3} + 2848z^{4}}{4z^{2}(-1 + 2z + 68z^{2} + 120z^{3})^{2}}$\\
\\
12B \> $\frac{1 - 4z + 54z^{2} - 388z^{3} + 769z^{4}}{4z^{2}(1 - 2z - 23z^{2} + 60z^{3})^{2}}$\\
\\
12C \> $\frac{1 + 32z^{2} + 336z^{4} + 10368z^{6}}{4z^{2}(1 - 40z^{2} + 144z^{4})^{2}}$\\
\\
12D \> $\frac{1 + 304z^{3} + 4096z^{6}}{4z^{2}(1 - 64z^{3})^{2}}$\\
\\
12E \> $\frac{1 + 720z^{4} + 1152z^{6}}{4z^{2}(1 + 8z^{2} - 48z^{4})^{2}}$\\
\\
12F \> $\frac{1 + 28z^{2} + 1350z^{4} + 28z^{6} + z^{8}}{4z^{2}(1 - 34z^{2} + z^{4})^{2}}$\\
\\
12G \> $\frac{1 - 4z^{2} - 474z^{4} - 324z^{6} + 6561z^{8}}{4z^{2}(1 + 14z^{2} + 81z^{4})^{2}}$\\
\\
12H \> $\frac{1 + 6z + 87z^{2} + 1908z^{3} + 14751z^{4} + 47142z^{5} + 54441z^{6}}{4z^{2}(-1 - 3z + 73z^{2} + 471z^{3} + 780z^{4})^{2}}$\\
\\
12I \> $\frac{(1 + 3z^{2})^{4}}{4z^{2}(1 - 10z^{2} + 9z^{4})^{2}}$\\
\\
12J \> $\frac{1 - 1600z^{6} + 4096z^{12}}{4z^{2}(1 + 64z^{6})^{2}}$\\
\\
13A \> $\frac{1 + 6z + 73z^{2} + 1512z^{3} + 11440z^{4} + 36048z^{5} + 41472z^{6}}{4z^{2}(-1 - 3z + 64z^{2} + 384z^{3} + 576z^{4})^{2}}$\\
\\
13B \> $\frac{1 + 2z + z^{2} + 74z^{3} + 224z^{4} + 1302z^{5} + 1345z^{6} + 3342z^{7} - 447z^{8}}{4z^{2}(1 + z + 8z^{2} - 9z^{3} - 3z^{4} - 70z^{5})^{2}}$\\
\\
14A \> $\frac{1 + 4z + 62z^{2} + 1036z^{3} + 7009z^{4} + 18888z^{5} + 17688z^{6}}{4z^{2}(-1 - 2z + 59z^{2} + 256z^{3} + 276z^{4})^{2}}$\\
\\
14B \> $\frac{1 - 2z + 15z^{2} - 180z^{3} + 1063z^{4} - 2802z^{5} + 2385z^{6} + 1544z^{7} - 1832z^{8}}{4z^{2}(-1 + z + 17z^{2} - 49z^{3} + 68z^{4} - 60z^{5})^{2}}$\\
\\
14C \> $\frac{1 + 12z + 106z^{2} + 1696z^{3} + 17779z^{4} + 98208z^{5} + 293690z^{6} + 454332z^{7} + 285777z^{8}}{4z^{2}(-1 - 6z + 45z^{2} + 490z^{3} + 1503z^{4} + 1564z^{5})^{2}}$\\
\\
15A \> $\frac{1 + 4z + 46z^{2} + 1040z^{3} + 6841z^{4} + 15628z^{5} + 10540z^{6}}{4z^{2}(-1 - 2z + 43z^{2} + 244z^{3} + 416z^{4})^{2}}$\\
\\
15B \> $\frac{1 - 2z^{2} + 152z^{3} - 357z^{4} + 488z^{5} + 1078z^{6} + 120z^{7} - 1175z^{8}}{4z^{2}(1 + 7z^{2} - 32z^{3} + 45z^{4} - 50z^{5})^{2}}$\\
\\
15C \> $\frac{1 + 10z + 83z^{2} + 1282z^{3} + 12073z^{4} + 59408z^{5} + 158768z^{6} + 220760z^{7} + 125800z^{8}}{4z^{2}(-1 - 5z + 43z^{2} + 377z^{3} + 965z^{4} + 825z^{5})^{2}}$\\
\\
15D \> $\frac{1 - 64z^{3} - 1750z^{6} - 8000z^{9} + 15625z^{12}}{4z^{2}(1 + 22z^{3} + 125z^{6})^{2}}$\\
\\
16A \> $\frac{1 + 16z^{2} + 608z^{4} + 256z^{6} + 256z^{8}}{4z^{2}(1 - 24z^{2} + 16z^{4})^{2}}$\\
\\
16B \> $\frac{1 + 224z^{4} + 256z^{8}}{4z^{2}(1 - 16z^{4})^{2}}$\\
\\
16C \> $\frac{1 + 4z + 44z^{2} + 800z^{3} + 4912z^{4} + 12864z^{5} + 12864z^{6}}{4z^{2}(-1 - 2z + 44z^{2} + 200z^{3} + 224z^{4})^{2}}$\\
\\
17A \> $\frac{1 + 8z + 58z^{2} + 864z^{3} + 7265z^{4} + 30984z^{5} + 70840z^{6} + 84016z^{7} + 41200z^{8}}{4z^{2}(-1 - 4z + 35z^{2} + 254z^{3} + 564z^{4} + 424z^{5})^{2}}$\\
\\
18A \> $\frac{1 - 4z^{2} + 40z^{3} - 156z^{4} + 896z^{5} - 176z^{6} + 3936z^{7} + 2024z^{8}}{4z^{2}(1 + 14z^{2} - 7z^{3} + 48z^{4} - 56z^{5})^{2}}$\\
\\
18B \> $\frac{1 + 2z + 35z^{2} + 468z^{3} + 2563z^{4} + 4978z^{5} + 3617z^{6}}{4z^{2}(1 + z - 39z^{2} - 109z^{3} - 70z^{4})^{2}}$\\
\\
18C \> $\frac{1 - 4z + 16z^{2} - 124z^{3} + 406z^{4} - 940z^{5} + 3520z^{6} - 8500z^{7} + 7225z^{8}}{4z^{2}(1 - 2z - 18z^{2} + 52z^{3} + z^{4} - 66z^{5})^{2}}$\\
\\
18D \> $\frac{1 + 40z^{3} + 384z^{6} - 320z^{9} + 64z^{12}}{4z^{2}(-1 + 7z^{3} + 8z^{6})^{2}}$\\
\\
18E \> $\frac{1 + 8z + 52z^{2} + 776z^{3} + 6532z^{4} + 27392z^{5} + 60688z^{6} + 68192z^{7} + 30472z^{8}}{4z^{2}(-1 - 4z + 30z^{2} + 227z^{3} + 524z^{4} + 420z^{5})^{2}}$\\
\\
19A \> $\frac{1 + 6z + 41z^{2} + 580z^{3} + 4120z^{4} + 14908z^{5} + 29212z^{6} + 29064z^{7} + 11340z^{8}}{4z^{2}(-1 - 3z + 32z^{2} + 172z^{3} + 276z^{4} + 144z^{5})^{2}}$\\
\\
20A \> $\frac{1 + 32z^{2} + 304z^{3} + 1088z^{4} + 2304z^{5} + 3072z^{6}}{4z^{2}(-1 + 32z^{2} + 64z^{3})^{2}}$\\
\\
20B \> $\frac{1 + 8z^{2} + 848z^{4} - 704z^{6} + 3072z^{8}}{4z^{2}(-1 + 12z^{2} + 64z^{4})^{2}}$\\
\\
20C \> $\frac{1 - 4z + 16z^{2} - 124z^{3} + 630z^{4} - 1292z^{5} + 1408z^{6} - 1172z^{7} + 793z^{8}}{4z^{2}(1 - 2z - 2z^{2} + 20z^{3} - 47z^{4} + 30z^{5})^{2}}$\\
\\
20D \> $\frac{1 - 64z^{4} - 576z^{6} - 1024z^{8} + 4096z^{12}}{4z^{2}(1 + 16z^{2} + 64z^{4} + 64z^{6})^{2}}$\\
\\
20E \> $\frac{1 + 10z^{2} + 351z^{4} + 44z^{6} + 351z^{8} + 10z^{10} + z^{12}}{4z^{2}(-1 + 19z^{2} - 19z^{4} + z^{6})^{2}}$\\
\\
20F \> $\frac{1 + 10z + 65z^{2} + 720z^{3} + 6350z^{4} + 32700z^{5} + 101870z^{6} + 198000z^{7} + 237105z^{8} + 161530z^{9} + 48337z^{10}}{4z^{2}(-1 - 5z + 20z^{2} + 210z^{3} + 627z^{4} + 835z^{5} + 426z^{6})^{2}}$\\
\\
21A \> $\frac{1 + 4z + 30z^{2} + 352z^{3} + 1705z^{4} + 5916z^{5} + 11340z^{6}}{4z^{2}(1 + 2z - 35z^{2} - 108z^{3})^{2}}$\\
\\
21B \> $\frac{1 + 2z - 5z^{2} - 54z^{3} + 9z^{4} + 576z^{5} + 576z^{6} + 792z^{7} + 504z^{8}}{4z^{2}(1 + z + 5z^{2} + 11z^{3} + 3z^{4} - 21z^{5})^{2}}$\\
\\
21C \> $\frac{1 + 110z^{3} + 1729z^{6} - 3672z^{9}}{4z^{2}(-1 + 26z^{3} + 27z^{6})^{2}}$\\
\\
21D \> $\frac{1 + 10z + 63z^{2} + 628z^{3} + 5257z^{4} + 26478z^{5} + 81681z^{6} + 158364z^{7} + 190143z^{8} + 130338z^{9} + 39321z^{10}}{4z^{2}(-1 - 5z + 21z^{2} + 197z^{3} + 511z^{4} + 561z^{5} + 222z^{6})^{2}}$\\
\\
22A \> $\frac{1 + 4z + 28z^{2} + 332z^{3} + 2144z^{4} + 6480z^{5} + 9548z^{6} + 5712z^{7} + 624z^{8}}{4z^{2}(-1 - 2z + 28z^{2} + 100z^{3} + 112z^{4} + 48z^{5})^{2}}$\\
\\
22B \> $\frac{1 - 4z + 12z^{2} - 100z^{3} + 736z^{4} - 2816z^{5} + 6636z^{6} - 9872z^{7} + 8656z^{8} - 3200z^{9} - 2240z^{10} + 3840z^{11} - 1280z^{12}}{4z^{2}(1 - 2z - 4z^{2} + 28z^{3} - 80z^{4} + 128z^{5} - 128z^{6} + 64z^{7})^{2}}$\\
\\
23A \> $\frac{1 + 8z + 44z^{2} + 462z^{3} + 3686z^{4} + 16546z^{5} + 44173z^{6} + 72694z^{7} + 72287z^{8} + 39090z^{9} + 8277z^{10}}{4z^{2}(-1 - 4z + 18z^{2} + 142z^{3} + 351z^{4} + 394z^{5} + 175z^{6})^{2}}$\\
\\
24A \> $\frac{1 + 8z^{2} + 512z^{6} + 4096z^{8}}{4(z - 20z^{3} + 64z^{5})^{2}}$\\
\\
24B \> $\frac{1 + 4z + 20z^{2} + 368z^{3} + 2288z^{4} + 6400z^{5} + 10240z^{6} + 10240z^{7} + 5120z^{8}}{4z^{2}(-1 - 2z + 16z^{2} + 96z^{3} + 192z^{4} + 128z^{5})^{2}}$\\
\\
24C \> $\frac{1 - 4z + 6z^{2} + 68z^{3} - 393z^{4} + 408z^{5} + 1716z^{6} - 5688z^{7} + 8223z^{8} - 9620z^{9} + 10982z^{10} - 8460z^{11} + 2825z^{12}}{4z^{2}(1 - 2z + z^{2} - 18z^{3} + 47z^{4} - 14z^{5} - 49z^{6} + 34z^{7})^{2}}$\\
\\
24D \> $\frac{1 - 144z^{4} - 1280z^{6} - 2304z^{8} + 4096z^{12}}{4z^{2}(1 + 8z^{2} + 32z^{4} + 64z^{6})^{2}}$\\
\\
24E \> $\frac{1 + 1600z^{6} + 4096z^{12}}{4z^{2}(1 - 64z^{6})^{2}}$\\
\\
24F \> $\frac{1 + 316z^{4} + 1926z^{8} + 316z^{12} + z^{16}}{4z^{2}(1 - 34z^{4} + z^{8})^{2}}$\\
\\
24G \> $\frac{1 - 100z^{4} - 2970z^{8} - 8100z^{12} + 6561z^{16}}{4z^{2}(1 + 14z^{4} + 81z^{8})^{2}}$\\
\\
24H \> $\frac{1 + 6z^{2} + 399z^{4} - 236z^{6} + 399z^{8} + 6z^{10} + z^{12}}{4z^{2}(1 - 13z^{2} - 13z^{4} + z^{6})^{2}}$\\
\\
24I \> \tiny $(1 + 12z + 78z^{2} + 612z^{3} + 4911z^{4} + 27384z^{5}
+ 100676z^{6} + 249672z^{7} + 422799z^{8} + 482940z^{9} +$ \\ \>
\tiny $ 355470z^{10} + 151764z^{11} + 28513z^{12}) / (4z^{2}(-1 - 6z
+ 11z^{2} + 174z^{3} + 553z^{4} + 798z^{5} + 541z^{6} +
138z^{7})^{2})$\\
\\
24J \> $\frac{1 - 6784z^{12}}{4z^{2} + 512z^{14}}$\\
\\
25A \>
$\frac{1+4z+22z^{2}+268z^{3}+1505z^{4}+4464z^5+7456z^{6}+6168z^{7} +
2112z^{8}}{4z^{2}(-1-2z+23z^{2}+83z^{3}+78z^{4}+19z^{5})^{2}}$\\
\\
26A \>
$\frac{1+4z+22z^{2}+228z^{3}+1345z^{4}+3888z^{5}+6624z^{6}+7376z^{7}+3840z^{8}}{4z^{2}(-1-2z+23z^{2}+76z^{3}+64z^{4})^{2}}$\\
\\
26B \> \tiny
$(1+14z+101z^{2}+740z^{3}+5908z^{4}+36852z^{5}+161298z^{6}+500052z^{7}+1122882z^{8}+1851100z^{9}
+ 2241876z^{10} + 1963740z^{11}$ \\ \> \tiny $+1192717z^{12} +
455190z^{13} +
83097z^{14})/(4z^{2}(-1-7z-2z^{2}+142z^{3}+674z^{4}+1536z^{5}+1959z^{6}+1341z^{7}+382z^{8})^{2})$\\
\\
27A \>
$\frac{1+6z+27z^{2}+288z^{3}+2079z^{4}+7938z^{5}+17577z^{6}+23328z^{7}+17496z^{8}+5832z^{9}}{4z^{2}(-1-3z+15z^{2}+90z^{3}+189z^{4}+189z^{5}+81z^{6})^2}$\\
\\
27B \>
$\frac{1+6z+27z^{2}+288z^{3}+2079z^{4}+7938z^{5}+17577z^{6}+23328z^{7}+17496z^{8}+5832z^{9}}{4z^{2}(-1-3z+15z^{2}+90z^{3}+189z^{4}+189z^{5}+81z^{6})^{2}}$\\
\\
28A \>
$\frac{1+4z^2-250z^{4}+3820z^{6}-19583z^{8}+31752z^{10}}{4z^{2}(1-22z^{2}+121z^{4}-196z^{6})^{2}}$\\
\\
28B \>
$\frac{1+2z+15z^{2}+180z^{3}+1063z^{4}+2802z^{5}+2385z^{6}-1544z^{7}-1832z^{8}}{4z^{2}(1+z-17z^{2}-49z^{3}-68z^{4}-60z^{5})^{2}}$\\
\\
28C \>
$\frac{1+2z^{2}+63z^{4}+292z^{6}+2223z^{8}+7538z^{10}+5425z^{12}+840z^{14}}{4z^{2}(1+9z^{2}+15z^{4}+3z^{6}-28z^{8})^{2}}$\\
\\
28D \>
$\frac{1+4z^{2}+170z^{4}+80z^{6}-909z^{8}+80z^{10}+170z^{12}+4z^{14}+z^{16}}{4z^{2}(1-14z^{2}+19z^{4}-14z^{6}+z^{8})^{2}}$\\
\\
29A \>
$\frac{1+8z+38z^{2}+296z^{3}+2091z^{4}+9000z^{5}+24526z^{6}+45520z^{7}+59625z^{8}+53152z^{9}+28472z^{10}+7248z^{11}+624z^{12}}{4z^{2}(-1-4z+13z^{2}+100z^{3}+215z^{4}+198z^{5}+68z^{6}+8z^{7})^{2}}$\\
\\
30A \> \tiny
$(1-10z+43z^{2}-158z^{3}+289z^{4}+2432z^{5}-15664z^{6}+23264z^{7}+61972z^{8}-318416z^{9}+633016z^{10}$\\
\> \tiny $-708784z^{11} +
368536z^{12})/(4z^{2}(-1+5z+15z^{2}-113z^{3}+57z^{4}+551z^{5}-904z^{6}+264z^{7})^{2})$\\
\\
30B \>
$\frac{1+2z+15z^{2}+112z^{3}+311z^{4}+1386z^{5}+2613z^{6}-2340z^{7}+2700z^{8}}{4z^{2}(1+z-25z^{2}-37z^{3}+60z^{4})^{2}}$\\
\\
30C \>
$\frac{1-2z+5z^{2}-76z^{3}+362z^{4}-1244z^{5}+2710z^{6}-3368z^{7}+4169z^{8}-2658z^{9}+7077z^{10}-5580z^{11}+2700z^{12}}{4z^{2}(-1+z-2z^{2}-18z^{3}+39z^{4}-79z^{5}+60z^{6})^{2}}$\\
\\
30D \> \tiny
$(1+10z+53z^{2}+362z^{3}+2414z^{4}+10642z^{5}+31681z^{6}+71734z^{7}+135662z^{8}+209214z^{9}+231861z^{10}$\\
\> \tiny
$+155646z^{11}+46521z^{12})/(4z^{2}(1+5z-10z^{2}-117z^{3}-282z^{4}-221z^{5}+69z^{6}+126z^{7})^{2})$\\
\\
30E \>
$\frac{1+70z^3+291z^6+1004z^{9}+291z^{12}+70z^{15}+z^{18}}{4z^{2}(-1+19z^{3}-19z^{6}+z^{9})^{2}}$\\
\\
30F \>
$\frac{1+6z+27z^{2}+210z^{3}+1497z^{4}+6072z^{5}+14360z^{6}+20448z^{7}+17172z^{8}+7600z^{9}+1080z^{10}-240z^{11}-40z^{12}}{4z^{2}(-1-3z+15z^{2}+75z^{3}+137z^{4}+139z^{5}+80z^{6}+20z^{7})^{2}}$\\
\\
30G \> \tiny
$(1-10z+47z^{2}-178z^{3}+853z^{4}-4148z^{5}+14908z^{6}-36712z^{7}+62456z^{8}-79616z^{9}+99968z^{10}$\\
\> \tiny
$-159232z^{11}+249824z^{12}-293696z^{13}+238528z^{14}-132736z^{15}+54592z^{16}-22784z^{17}+12032z^{18}$\\
\> \tiny
$-5120z^{19}+1024z^{20})/(4z^{2}(1-5z+3z^{2}+33z^{3}-115z^{4}+9205z^{5}-230z^{6}+132z^{7}+24z^{8}$\\
\> \tiny $-80z^{9}+36z^{10})^{2})$\\
\\
31A \>
$\frac{1+8z+36z^{2}+246z^{3}+1622z^{4}+6930z^{5}+18981z^{6}+33086z^{7}+35047z^{8}+22122z^{9}+10189z^{10}+3864z^{11}+648z^{12}}{4z^{2}(-1-4z+14z^{2}+94z^{3}+159z^{4}+98z^{5}+27z^{6})^{2}}$\\
\\
32A \>
$\frac{1+8z+36z^{2}+272z^{3}+2044z^{4}+9600z^{5}+28640z^{6}+56576z^{7}+74096z^{8}+60800z^{9}+25664z^{10}+768z^{11}-2496z^{12}}{4z^{2}(-1-4z+6z^{2}+76z^{3}+228z^{4}+352z^{5}+296z^{6}+112z^{7})^{2}}$\\
\\
33A \> \tiny
$(1-2z+5z^{2}+30z^{3}-214z^{4}+790z^{5}-1843z^{6}+3966z^{7}-5586z^{8}+9038z^{9}-9919z^{10}$\\
\> \tiny
$+13318z^{11}-11295z^{12}+10912z^{13}-5368z^{14}+4312z^{15}-1584z^{16})/(4z^{2}(-1+z-10z^{2}+18z^{3}$\\
\> \tiny
$-55z^{4}+87z^{5}-122z^{6}+126z^{7}-77z^{8}+33z^{9})^{2})$\\
\\
33B \>
$\frac{1+4z+14z^{2}+188z^{3}+1081z^{4}+2756z^{5}+4976z^{6}+8092z^{7}+9200z^{8}+5104z^{9}+1452z^{10}}{4z^{2}(-1-2z+11z^{2}+56z^{3}+96z^{4}+44z^{5})^{2}}$\\
\\
34A \>
$\frac{1+2z+11z^{2}+100z^{3}+411z^{4}+898z^{5}+1417z^{6}+752z^{7}-664z^{8}-912z^{9}+432z^{10}}{4z^{2}(1+z-19z^{2}-33z^{3}+14z^{4}+20z^{5}-8z^{6})^{2}}$\\
\\
35A \>
$\frac{1+4z+14z^{2}+180z^{3}+1273z^{4}+4780z^{5}+11632z^{6}+15964z^{7}+12112z^{8}-3360z^{9}-6900z^{10}}{4z^{2}(-1-2z+3z^{2}+44z^{3}+136z^{4}+220z^{5}+200z^{6})^{2}}$\\
\\
35B \> \tiny $(1 + 10z + 53z^{2} + 314z^{3} + 2130z^{4} + 10970z^{5}
+ 39221z^{6} + 100130z^{7} + 188542z^{8} + 267514z^{9}$\\

\> \tiny $+290017z^{10} + 243226z^{11} +
160073z^{12}+83320z^{13}+33272z^{14}+9048z^{15}+1224z^{16})/(4z^{2}(-1-5z$\\
\> \tiny
$+2z^{2}+74z^{3}+257z^{4}+465z^{5}+510z^{6}+342z^{7}+127z^{8}+19z^{9})^{2})$\\
\\
36A \>
$\frac{1+4z+16z^{2}+124z^{3}+406z^{4}+940z^{5}+3520z^{6}+8500z^{7}+7225z^{8}}{4z^{2}(1+2z-18z^{2}-52z^{3}+z^{4}+66z^{5})^{2}}$\\
\\
36B \>
$\frac{1-40z^{3}+384z^{6}+320z^{9}+64z^{12}}{4z^{2}(1+7z^{3}-8z^{6})^{2}}$\\
\\
36C \>
$\frac{1+2z^{2}+279z^{4}-612z^{6}+2511z^{8}+162z^{10}+729z^{12}}{4z^{2}(1-7z^{2}-21z^{4}+27z^{6})^{2}}$\\
\\
36D \> \tiny
$(1+12z+72z^{2}+40z^{3}+2580z^{4}+13536z^{5}+51696z^{6}+144864z^{7}+308376z^{8}+514080z^{9}$\\
\> \tiny
$+682560z^{10}+720576z^{11}+591840z^{12}+362880z^{13}+155520z^{14}+41472z^{15}+5184z^{16})/(4z^{2}(-1-6z$\\
\> \tiny
$-2z^{2}+81z^{3}+318z^{4}+588z^{5}+600z^{6}+324z^{7}+72z^{8})^{2})$\\
\\
38A \>
$\frac{1+4z+14z^{2}+112z^{3}+769z^{4}+2440z^{5}+4532z^{6}+5260z^{7}+3976z^{8}+3328z^{9}+3308z^{10}+2464z^{11}+768z^{12}}{4z^{2}(-1-2z+11z^{2}+40z^{3}+68z^{4}+60z^{5}+32z^{6})^{2}}$\\
\\
39A \>
$\frac{1+8z^{2}+42z^{3}+792z^{5}+161z^{6}-4104z^{7}+5608z^{8}-3048z^{9}+648z^{10}}{4z^{2}(1-20z^{2}-6z^{3}+64z^{4}-48z^{5}+9z^{6})^{2}}$\\
\\
39B \>
$\frac{1+36z^{3}+1006z^{6}+828z^{9}-2519z^{12}-4104z^{15}}{4z6{2}(-1+9z^{3}+37z^{6}+27z^{9})^{2}}$\\
\\
39C \> \tiny
$(1+10z+51z^{2}+260z^{3}+1567z^{4}+7836z^{5}+27877z^{6}+70044z^{7}+126450z^{8}+166144z^{9}+158719z^{10}$\\
\> \tiny
$+106092z^{11}+40861z^{12}-4242z^{13}-19146z^{14}-14760z^{15} -
6435z^{16}-1620z^{17} - 180z^{18})/(4z^{2}(-1-5z$\\
\> \tiny $+3z^{2}+69z^{3}+200z^{4}+308z^{5}+310z^{6}+216z^{7}+103z^{8}+33z^{9}+6z^{10})^{2})$\\
\\
40A \>
$\frac{1+104z^{4}+3344z^{8}-7424z^{12}}{4z^{2}(-1+12z^{4}+64z^{8})^{2}}$\\
\\
40B \>
$\frac{1-64z^{4}+576z^{6}-1024z^{8}+4096z^{12}}{4z^{2}(1-16z^{2}+64z^{4}-64z^{6})^{2}}$\\
\\
40C \>
$\frac{1+2z^{2}+173z^{4}+280z^{6}+82z^{8}-1076z^{10}+82z^{12}+280z^{14}+173z^{16}+2z^{18}+z^{20}}{4z^{2}(1-7z^{2}-10z^{4}-10z^{6}-7z^{8}+z^{10})^{2}}$\\
\\
41A \> \tiny
$(1+8z+32z^{2}+168z^{3}+992z^{4}+3880z^{5}+9914z^{6}+18616z^{7}+27584z^{8}
+ 28440z^{9} + 11872z^{10}$\\ \> \tiny
$-8504z^{11}-8839z^{12}+848z^{13}+1184z^{14}-624z^{15}+192z^{16})/(4z^{2}(-1-4z+8z^{2}+66z^{3}+120z^{4}$\\
\> \tiny $+56z^{5} - 53z^{6}-36z^{7}+16z^{8})^{2})$\\
\\
42A \>
$\frac{1+2z+7z^{2}+88z^{3}+343z^{4}+122z^{5}+781z^{6}+3108z^{7}+1020z^{8}-288z^{9}+1728z^{10}}{4z^{2}(1+z-13z^{2}-29z^{3}-8z^{4}+48z^{5})^{2}}$\\
\\
42B \> \tiny
$(1-10z+43z^{2}-106z^{3}+133z^{4}-260z^{5}+4730z^{6}-31616z^{7}+97282z^{8}-131528z^{9}-39382z^{10}$\\
\> \tiny
$+407916z^{11}-576411z^{12}+277146z^{13}+117675z^{14}-183726z^{15}+58257z^{16})/(4z^{2}(1-5z+z^{2}+32z^{3}$\\
\> \tiny
$-46z^{4}-14z^{5}-79z^{6}+503z^{7}-711z^{8}+330z^{9})^{2})$\\
\\
42C \>
$\frac{1-26z^{3}-519z^{6}-4440z^{9}-18368z^{12}-35520z^{15}-33216z^{18}-13312z^{21}+4096z^{24}}{4z^{2}(1+14z^{3}+61z^{6}+112z^{9}+64z^{12})^{2}}$\\
\\
42D \> \tiny
$(1+8z+34z^{2}+212z^{3}+1411z^{4}+6508z^{5}+21074z^{6}+50668z^{7}+92701z^{8}+128728z^{9}+132692z^{10}$\\
\> \tiny
$+97272z^{11}+46788z^{12}+12096z^{13}+288z^{14}-576z^{15}-72z^{16})/(4z^{2}(-1-4z-z^{2}+43z^{3}+163z^{4}$\\
\> \tiny $+311z^{5}+364z^{6}+271z^{7}+120z^{8}+24z^{9})^{2})$\\
\\
44A \>
$\frac{1+4z+12z^{2}+100z^{3}+736z^{4}+2816z^{5}+6636z^{6}+9872z^{7}+8656z^{8}+3200z^{9}-2240z^{10}-3840z^{11}-1280z^{12}}{4z^{2}(-1-2z+4z^{2}+28z^{3}+80z^{4}+128z^{5}+128z^{6}+64z^{7})^{2}}$\\
\\
45A \>
$\frac{1+4z+10z^{2}+64z^{3}+363z^{4}+1416z^{5}+2850z^{6}+1524z^{7}-3543z^{8}-7232z^{9}-3776z^{10}-8z^{11}+592z^{12}}{4z^{2}(-1-2z+13z^{2}+33z^{3}+9z^{4}+13z^{5}+34z^{6}+11z^{7})^{2}}$\\
\\
46A \> \tiny
$(1-6z+23z^{2}-94z^{3}+431z^{4}-1834z^{5}+6210z^{6}-16564z^{7}+35555z^{8}-63012z^{9}+96486z^{10}$\\
\> \tiny
$-132532z^{11}+171542z^{12}-206124z^{13}+221548z^{14}-189264z^{15}+106711z^{16}+1034z^{17}-77957z^{18}$\\
\> \tiny
$+96226z^{19}-66669z^{20}+28278z^{21}-6694z^{22}+563z^{23}+237z^{24})/(4z^{2}(1-3z+7z^{2}+z^{3}-32z^{4}$\\
\> \tiny
$+110z^{5}-197z^{6}+247z^{7}-198z^{8}+90z^{9}-14z^{10}-4z^{11}-15z^{12}+7z^{13})^{2})$\\
\\
46C \> \tiny
$(1+4z+10z^{2}+62z^{3}+359z^{4}+1174z^{5}+1837z^{6}+484z^{7}-609z^{8}+2520z^{9}+1612z^{10}
- 3378z^{11}$\\ \> \tiny
$-863z^{12}+414z^{13}+45z^{14})/(4z^{2}(-1-2z+13z^{2}+34z^{3}+9z^{4}-10z^{5}+2z^{6}-12z^{7}+7z^{8})^{2})$\\
\\
47A \> \tiny $(1+8z+36z^{2}+188z^{3} + 1158z^{4} + 5570z^{5} +
19584z^{6} + 51414z^{7} + 103965z^{8} + 165234z^{9} +
209397z^{10}$\\ \> \tiny $+213722z^{11} + 177346z^{12} +
120520z^{13} + 66397z^{14} + 27640z^{15} + 6488z^{16} - 648z^{17} -
708z^{18})/(4z^{2}(-1-4z$\\ \> \tiny $-2z^{2} + 32z^{3} + 135z^{4} +
294z^{5} + 424z^{6} + 410z^{7}+268z^{8}+100z^{9} + 19z^{10})^{2})$\\
\\
48A \>
$\frac{1+64z^{4}+640z^{6}-2592z^{8}+2560z^{10}+1024z^{12}+256z^{16}}{4z^{2}(1-8z^{2}+8z^{4}-32z^{6}+16z^{8})^{2}}$\\
\\
50A \>
$\frac{1+4z+10z^{2}+60z^{3}+245z^{4}+496z^{5}+1124z^{6}+2840z^{7}+2380z^{8}
-
3040z^{9}-4096z^{10}+1536z^{11}+2280z^{12}}{4z^{2}(1+2z-13z^{2}-35z^{3}+53z^{5}+16z^{6}-24z^{7})^{2}}$\\
\\
51A \> \tiny
$(1+4z+10z^{2}+92z^{3}+499z^{4}+1308z^{5}+2598z^{6}+4388z^{7}+5905z^{8}+5212z^{9}+2024z^{10}$\\
\> \tiny $-2628z^{11} -
4288z^{12}-1536z^{13}-116z^{14})/(4z^{2}(-1-2z+5z^{2}+30z^{3}+55z^{4}+36z^{5}+8z^{6}$\\
\> \tiny $-4z^{7} + 8z^{8})^{2})$\\
\\
52A \>$
\frac{1+4z^{2}+214z^{4}+1644z^{6}+5089z^{8}+5832z^{10}+336z^{12}+320z^{14}+3072z^{16}}{4z^{2}(-1-2z^{2}+23z^{4}+76z^{6}+64z^{8})^{2}}$\\
\\
52B \> \tiny $(1+2z^{2} + 65z^{4}+108z^{6}-1224z^{8} + 3308z^{10} -
4890z^{12} + 6036z^{14} - 4890z^{16} + 3308z^{18}$\\ \> \tiny
$-1224z^{20} + 108z^{22} + 65z^{24} - 2z^{26} + z^{28})/(4z^{2}(-1 +
9z^{2} - 16z^{4} + 26z^{6}-26z^{8} + 16z^{10}$\\ \> \tiny $-9z^{12}
+ z^{14})^{2})$\\
\\
54A \> \tiny $(1+4z+12z^{2} + 64z^{3} + 478z^{4} + 1788z^{5} +
4540z^{6} + 8272z^{7} + 11169z^{8} + 10880z^{9} + 6404z^{10}$\\ \>
\tiny $-696z^{11} - 6076z^{12} - 8128z^{13} - 5712z^{14} -
2848z^{15} - 712z^{16})/(4z^{2}(-1-2z+4z^{2}+19z^{3} + 55z^{4}$\\ \>
\tiny $+88z^{5} + 106z^{6} + 91z^{7} + 52z^{8} + 20z^{9})^{2})$\\
\\
55A \> \tiny $(1+4z+6z^{2} + 46z^{3} + 105z^{4}-160z^{5} + 201z^{6}
+ 3114z^{7} + 3154z^{8} - 6500z^{9} - 15911z^{10}$\\ \> \tiny
$-8320z^{11} + 12400z^{12} + 17960z^{13} + 7640z^{14})/(4z^{2}(-1 -
2z+15z^{2} + 38z^{3} - 36z^{4} - 126z^{5} + 7z^{6}$\\ \> \tiny
$+110z^{7} + 35z^{8})^{2})$\\
\\
56A \> \tiny $(1 + 6z+19z^{2}+112z^{3} + 569z^{4} + 1706z^{5} +
3167z^{6} + 4900z^{7} + 9735z^{8} + 20058z^{9} + 29657z^{10}$\\ \>
\tiny $+29448z^{11} + 22963z^{12} + 15798z^{13} + 8693z^{14} +
3108z^{15} + 588z^{16})/(4z^{2}(1+3z-3z^{2} - 37z^{3} -77z^{4}$\\ \>
\tiny $-51z^{5} + 51z^{6}+85z^{7} + 28z^{8})^{2})$\\
\\
57A \> $\frac{1+30z^{3}-z^{6}-728z^{9} + 6691z^{12} - 17630z^{15} +
13753z^{18} -
216z^{21}}{4z^{2}(1-12z^{3}+22z^{6}+8z^{9}-27z^{12})^{2}}$\\
\\
59A \> \tiny $(1+8z+32z^{2} + 126z^{3} + 664z^{4} + 3064z^{5} +
10303z^{6} + 25672z^{7} + 48984z^{8} + 73000z^{9} + 85432z^{10}$\\
\> \tiny $+75264z^{11} + 39939z^{12} - 9016z^{13} - 53200z^{14} -
78750z^{15} - 80728z^{16} - 65192z^{17} - 42791z^{18} -
22808z^{19}$\\ \> \tiny
$-9496z^{20}-2936z^{21}-536z^{22})/(4z^{2}(-1 - 4z + 28z^{3} +
84z^{4} + 152z^{5} + 202z^{6} + 212z^{7} + 176z^{8} + 120z^{9}$\\ \>
\tiny $+68z^{10} + 24z^{11} + 11z^{12})^{2})$\\
\\
60A \> $\frac{1-6z^{2}-169z^{4} + 2612z^{6}-15153z^{8} + 43978z^{10}
- 71655z^{12} + 68040z^{14}}{4z^{2}(1-19z^{2}+119z^{4}-281z^{6} +
180z^{8})^{2}}$\\
\\
60B \>
$\frac{1+2z+5z^{2}+76z^{3}+362z^{4}+1244z^{5}+2710z^{6}+3368z^{7}+4169z^{8}+2658z^{9}
+ 7077z^{10} + 5580z^{11} +
2700z^{12}}{4z^{2}(1+z+2z^{2}-18z^{3}-39z^{4}-79z^{5}-60z^{6})^{2}}$\\
\\
60C \> \tiny $(1+10z+47z^{2}+178z^{3}+853z^{4}+4148z^{5} +
14908z^{6} + 36712z^{7} + 62456z^{8} + 79616z^{9} + 99968z^{10}$\\
\> \tiny $+159232z^{11} + 249824z^{12} + 293696z^{13} + 238528z^{14}
+ 132736z^{15} + 54592z^{16} + 22784z^{17} + 12032z^{18} +
5120z^{19}$\\ \> \tiny
$+1024z^{20})/(4z^{2}(1+5z+3z^{2}-33z^{3}-115z^{4}-205z^{5}-230z^{6}-132z^{7}+24z^{8}
+ 80z^{9} + 32z^{10})^{2})$\\
\\
60D \> \tiny $(1-6z+19z^{2}-10z^{3}-167z^{4} +
1120z^{5}-4576z^{6}+13408z^{7} - 30476z^{8} + 56192z^{9} -
81608z^{10} + 80272z^{11}$\\ \> \tiny $-21144z^{12} - 93056z^{13} +
211328z^{14} - 274432z^{15} + 261152z^{16} - 193152z^{17} +
111040z^{18} - 48000z^{19} + 14545z^{20}$\\ \> \tiny $-2816z^{21} +
256z^{22})/(4z^{2}(1-3z+13z^{2}-41z^{3}+95z^{4}-197z^{5}+344z^{6}-452z^{7}+420z^{8}
- 260z^{9}$\\ \> \tiny$+96z^{10} - 16z^{11})^{2})$\\
\\
60E \>
$\frac{1-2z^{2}+61z^{4}-294z^{6}+122z^{8}+3446z^{10}-5003z^{12}+3446z^{14}
+ 122z^{16} - 294z^{18} + 61z^{20}-2z^{22} +
z^{24}}{4z^{2}(1-9z^{2}+14z^{4}+3z^{6}+14z^{8}-9z^{10}+z^{12})^{2}}$\\
\\
60F \>
$\frac{1+394z^{6}-33z^{12}+6188z^{18}-33z^{24}+394z^{30}+z^{36}}{4z^{2}(-1+19z^{6}-19z^{12}+z^{18})^{2}}$\\
\\
62A \> \tiny $(1+4z+10z^{2}+54z^{3}+343z^{4} + 998z^{5} + 1821z^{6}
+ 2084z^{7}+2271z^{8} + 4360z^{9} + 8436z^{10} + 9862z^{11}$\\ \>
\tiny $+4081z^{12}-4002z^{13}-7211z^{14}-3184z^{15}+304z^{16} +
1536z^{17} + 336z^{18})/(4z^{2}(-1-2z+5z^{2}+22z^{3}$\\ \> \tiny
$+41z^{4}+30z^{5}+2z^{6}-32z^{7}-21z^{8}-4z^{9}+12z^{10})^{2})$\\
\\
66A \> \tiny
$(1+2z-z^{2}+4z^{3}-57z^{4}+334z^{5}+1577z^{6}-3072z^{7}-7216z^{8} +
12476z^{9} + 3836z^{10} - 17208z^{11}$\\ \> \tiny $+20748z^{12} -
175628z^{13} +
6912z^{14})/(4z^{2}(1+z-17z^{2}-13z^{3}+80z^{4}-148z^{6} +
96z^{7})^{2})$\\
\\
66B \> \tiny
$(1+10z+47z^{2}+178z^{3}+853z^{4}+4148z^{5}+14908z^{6}+36712z^{7}+62456z^{8}
+ 79616z^{9} + 99968z^{10}$\\ \> \tiny $+159232z^{11}+249824z^{12} +
293696z^{13} + 238528z^{14} + 132736z^{15} + 54592z^{16} +
22784z^{17} + 12032z^{!8} + 5120z^{19}$\\ \> \tiny
$+1024z^{20})/(4z^{2}(1+5z+3z^{2}-33z^{3}-115z^{4}-205z^{5}-230z^{6}-132z^{7}+
24z^{8} + 80z^{9} + 32z^{10})^{2})$\\
\\
68A \>
$\frac{1-6z^{2}+71z^{4}-940z^{6}+6591z^{8}-24118z^{10}+47017z^{12}-47384z^{14}+21008z^{16}-2240z^{18}+3072z^{20}}{4z^{2}(1-11z^{2}+31z^{4}-z^{6}
- 84z^{8} + 64z^{10})^{2}}$\\
\\
69A \> \tiny
$(1+4z+6z^{2}+42z^{3}+217z^{4}+398z^{5}+719z^{6}+2344z^{7} +
1952z^{8} - 5504z^{9} - 7386z^{10} + 4220z^{11}$\\ \> \tiny
$+5291z^{12} - 6276z^{13} - 4794z^{14} + 3582z^{15} + 1233z^{16} -
2394z^{17} + 117z^{18})/(4z^{2}(-1-2z + 7z^{2} + 24z^{3}$\\ \> \tiny
$+16z^{4} - 22z^{5} - 10z^{6}+44z^{7} + 4z^{8} - 24z^{9} +
15z^{10})^{2})$\\
\\
70A \> \tiny $(1+ 2z + 3z^{2} + 8z^{3} + 179z^{4} + 806z^{5} +
965z^{6} + 676z^{7} + 2004z^{8} - 420z^{9} - 4468z^{10} -
2392z^{11}$\\ \> \tiny $-3652z^{12} - 8976z^{13} -
528z^{14})/(4z^{2}(-1 - z + 7z^{2} + 5z^{3} + 10z^{4} + 48z^{5} +
20z^{6} + 32z^{7} + 48z^{8})^{2})$\\
\\
70B \> \tiny $(1 - 8z + 32z^{2} - 116z^{3} + 518z^{4} - 2076z^{5} +
5912z^{6} - 10196z^{7} + 3463z^{8} + 35596z^{9} - 112764z^{10}$\\ \>
\tiny $+168492z^{11} - 58834z^{12} - 356588z^{13} + 1048384z^{14} -
1754936z^{15} + 2145569z^{16} - 2062200z^{17} + 1601820z^{18}$\\ \>
\tiny $-1007240z^{19} + 502972z^{20} - 191840z^{21} + 52400z^{22} -
9120z^{23} + 760z^{24})/(4z^{2}(1 - 4z + 8z^{2} + z^{3} - 49z^{4}$\\
\> \tiny $+134z^{5} - 176z^{6} + 35z^{7} + 325z^{8} - 670z^{9} +
694z^{10} - 419z^{11} + 140z^{12} - 20z^{13})^{2})$\\
\\
71A \> \tiny $(1 + 8z + 28z^{2} + 90z^{3} + 410z^{4} + 1458z^{5} +
3073z^{6} + 3798z^{7} + 4021z^{8} + 8108z^{9} + 13349z^{10} +
88z^{11}$\\ \> \tiny $-40168z^{12} - 63932z^{13} - 13304z^{14} +
79626z^{15} + 94791z^{16} - 5466z^{17} - 91159z^{18} - 52578z^{19} +
31090z^{20}$\\ \> \tiny $+41224z^{21} + 1849z^{22} - 12328z^{23} -
3112z^{24} + 1320z^{25} + 364z^{26})/(4z^{2}(1 + 4z - 2z^{2} -
38z^{3} - 77z^{4}$\\ \> \tiny $-26z^{5} + 111z^{6} + 148z^{7} +
z^{8} - 122z^{9} - 70z^{10} + 30z^{11} + 40z^{12} + 4z^{13} -
11z^{14})^{2})$\\
\\
78A \> \tiny $(1+2z + z^{2} + 34z^{3} + 132z^{4} - 154z^{5} -
591z^{6} + 398z^{7} + 2909z^{8} + 3656z^{9} - 2400z^{10} -
12752z^{11}$\\ \> \tiny $-10660z^{12} + 5648z^{13} + 14424z^{14} +
5904z^{15} + 648z^{16})(4z^{2}(1 + z - 8z^{2} - 18z^{3}-2z^{4} +
48z^{5} + 45z^{6}$\\ \> \tiny $-7z^{7} - 48z^{8} - 36z^{9})^{2})$\\
\\
78B \> \tiny $(1 - 10z + 49z^{2} - 156z^{3} + 364z^{4} - 890z^{5} +
3626z^{6} - 16788z^{7} + 62992z^{8} - 185826z^{9} + 442843z^{10}$\\
\> \tiny $-871628z^{11} + 1438554z^{12} - 2020966z^{13} +
2462186z^{14} - 2646524z^{15} + 2513284z^{16} - 2051814z^{17} +
1362881z^{18}$\\ \> \tiny $-692340z^{19} + 283888z^{20} -
171710z^{21} + 188800z^{22} - 176220z^{23} + 111405z^{24} -
46800z^{25} + 12672z^{26} - 2016z^{27}$\\ \> \tiny
$+144z^{28})/(4z^{2}(1 - 5z + 12z^{2} - 18z^{3} + 20z^{4} - 4z^{5} -
75z^{6} + 267z^{7} - 567z^{8} + 901z^{9} - 1129z^{10}$\\ \> \tiny
$+1083z^{11} - 749z^{12} + 347z^{13} - 96z^{14} + 12z^{15})^{2})$\\
\\
84A \> $\frac{1+2z^{2} + 103z^{4} + 604z^{6} + 1087z^{8} -
2014z^{10} - 311z^{12} + 14784z^{14} + 12240z^{16} +
1152z^{18}}{4z^{2}(1 + z^{2} - 13z^{4} - 29z^{6} - 8z^{8} +
48z^{10})^{2}}$\\
\\
84B \> \tiny $(1 + 10z^{2} - 5z^{4} - 506z^{6} - 3191z^{8} -
12204z^{10} - 37490z^{12} - 82088z^{14} - 108062z^{16} -
82088z^{18}$\\ \> \tiny $-37490z^{20} - 12204z^{22} - 3191z^{24} -
506z^{26} - 5z^{28} + 10z^{30} + z^{32})/(4z^{2}(1 + 13z^{2} +
65z^{4} + 158z^{6}$\\ \> \tiny $+198z^{8} + 158z^{10} + 65z^{12} +
13z^{14} + z^{16})^{2})$\\
\\
84C \> $\frac{1+26z^{3} - 519z^{6} + 4440z^{9} - 18368z^{12} +
35520z^{15} - 33216z^{18} + 13312z^{21} + 4096z^{24}}{4z^{2}(1-
14z^{3} + 61z^{6} - 112z^{9} + 64z^{12})^{2}}$\\
\\
87A \> \tiny $(1 + 4z + 14z^{2} + 66z^{3} + 305z^{4} +
1226z^{5}+4151z^{6} + 10852z^{7} + 24748z^{8} + 44920z^{9} +
72122z^{10}$\\ \> \tiny $+91948z^{11} + 101547z^{12} + 80524z^{13}
+ 44982z^{14} - 6562z^{15} - 39367z^{16} - 55542z^{17} - 45531z^{18}
- 29628z^{19}$\\ \> \tiny $-14220z^{20} - 5040z^{21} - 2124z^{22}) /
(4z^{2}(-1 -2z - 5z^{2} + 4z^{3} + 24z^{4} + 78z^{5} + 146z^{6} +
200z^{7} + 240z^{8}$\\ \> \tiny $+188z^{9} + 151z^{10} + 60z^{11} +
24z^{12})^{2})$\\
\\
88A \> \tiny $(1 - 8z^{2} - 16z^{4} + 840z^{6} - 6464z^{8} +
25792z^{10} - 64128z^{12} +103168z^{14} - 103424z^{16} +
53760z^{18}$\\ \> \tiny$-4096z^{20} - 8192z^{22} +
4096z^{24})/(4z^{2}(1-12z^{2} + 56z^{4} - 148z^{6} + 224z^{8} -
192z^{10} + 64z^{12})^{2})$\\
\\
92A \> \tiny $(1 + 6z + 23z^{2} + 94z^{3} + 431z^{4} + 1834z^{5} +
6210z^{6} + 16564z^{7} + 35555z^{8} + 63012z^{9} +96486z^{10}$\\
\> \tiny $+132532z^{11} + 171542z^{12} + 206124z^{13} +
221548z^{14} + 189264z^{15}+106711z^{16} - 1034z^{17} - 77957z^{18}$
\\ \> \tiny $-96226z^{19} - 66669z^{20} - 28278z^{21} - 6694z^{22} -
564z^{23} + 237z^{24})/(4z^{2}(-1 - 3z - 7z^{2} + z^{3} + 32z^{4}$\\
\> \tiny $+110z^{5} + 197z^{6} + 247z^{7} + 198z^{8} + 90z^{9} +
14z^{10} - 4z^{11} + 15z^{12} + 7z^{13})^{2})$\\
\\
93A \> \tiny $(1 + 8z^{3} + 420z^{6} + 4614z^{9} + 21590z^{12} +
53058z^{15} + 71229z^{18} + 38414z^{21} - 40265z^{24} -
90822z^{27}$\\ \> \tiny $-53819z^{30} - 7560z^{33})/(4z^{2}(-1 -
4z^{3}+14z^{6} + 94z^{9} + 159z^{12} + 98z^{15} + 27z^{18})^{2})$\\
\\
94A \> \tiny $(1 + 4z + 2z^{2} - 4z^{3} + 71z^{4} + 486z^{5} +
844z^{6} - 1026z^{7} - 2089z^{8} + 5456z^{9} + 2155z^{10}$\\ \>
\tiny $-18040z^{11} + 9269z^{12} + 28398z^{13} - 45506z^{14} -
7458z^{15} + 68336z^{16} - 55696z^{17} - 26471z^{18}$\\ \> \tiny
$+69872z^{19} - 43972z^{20} - 9720z^{21} + 29632z^{22} - 22032z^{23}
+ 7216z^{24} - 992z^{25} - 232z^{26})/(4z^{2}(1 + 2z$\\ \> \tiny
$-9z^{2} -16z^{3} + 19z^{4} - 59z^{6} + 40z^{7} + 20z^{8} - 90z^{9}
+ 65z^{10} + 6z^{11} - 63z^{12} + 52z^{!3} - 20z^{14})^{2})$\\
\\
95A \> \tiny $(1 + 4z + 2z^{2} - 4z^{3} + 67z^{4} + 478z^{5} +
848z^{6} - 1034z^{7} - 3347z^{8} + 1542z^{9} + 7379z^{10} +
348z^{11}$\\ \> \tiny $-6511z^{12} + 7508z^{13} + 6564z^{14} -
19472z^{15} + 9653z^{16} + 10646z^{17} - 18870z^{18} + 6050z^{19} +
3245z^{20}$\\ \> \tiny $-930z^{21} + 165z^{22})/(4z^{2}(1 + 2z -
9z^{2} - 16z^{3} + 17z^{4} - 43z^{6} + 36z^{7} + 27z^{8} - 22z^{9} +
2z^{!0} + 20z^{11}$\\ \> \tiny $-15z^{12})^{2})$\\
\\
104A \> $\frac{1+4z^{4} + 982z^{8} + 7308z^{12} + 20065z^{16} +
13608z^{20}-24816z^{24} -
27904z^{28}}{4z^{2}(-1-2z^{4}+23z^{8}+76z^{12} + 64z^{16})^{2}}$\\
\\
105A \> \tiny $(1+2z - 3z^{2} + 24z^{3} + 24z^{4} - 496z^{5} - 566z^{6} +
3432z^{7} + 5926z^{8} - 10512z^{9} - 26620z^{10} + 9344z^{11}$\\ \>
\tiny $+57793z^{12} + 17382z^{13} - 57447z^{14} - 44064z^{15} +
11880z^{16} + 30600z^{17} + 18000z^{18})/(4z^{2}(1+z-10z^{2}$\\ \>
\tiny $-21z^{3} + 25z^{4} + 100z^{5} + 19z^{6} - 1273z^{7} -
122z^{8} + 105z^{9} + 75z^{10})^{2})$\\
\\
110A \> \tiny $(1 + 2z + 7z^{2} + 38z^{3} + 179z^{4} + 440z^{5} +
1002z^{6} + 44z^{7} + 1091z^{8} - 3578z^{9} + 10043z^{10} +
1742z^{11}$\\ \> \tiny $+29593z^{12} - 16792z^{13} + 2440z^{14} -
35328z^{15} +15300z^{16} + 18128z^{17} + 4568z^{18} - 14736z^{19} +
5016z^{20})/4z^{2}(1$\\ \> \tiny $+z + 3z^{2} - 11z^{3} - 22z^{4} -
42z^{5} - 21z^{6} + 57z^{7} + 31z^{8} + 91z^{9} - 112z^{10} +
24z^{11})^{2})$\\
\\
119A \> \tiny $(1 + 4z + 14z^{2} + 40z^{3} + 193z^{4} + 814z^{5} +
3046z^{6} + 9970z^{7} + 27348z^{8} + 65380z^{9} + 137087z^{10}$\\ \>
\tiny $+ 251940z^{11} + 411609z^{12} + 591020z^{13} + 748046z^{14} +
828256z^{15} + 778011z^{16} + 618670z^{17} + 394924z^{18}$\\ \>
\tiny $+179710z^{19} + 62284z^{20} + 36940z^{21} + 41961z^{22} +
43232z^{23} + 34104z^{24} + 12824z^{25} + 868z^{26})/4z^{2}(-1$\\ \>
\tiny $-2z - 5z^{2} -10z^{3} + 24z^{5} + 91z^{6} + 230z^{7} +
391z^{8} + 560z^{9} + 659z^{10} + 596z^{11} + 429z^{12} +
182z^{13}$\\ \> \tiny $+35z^{14})^{2})$\\
\\
25z \> $\frac{1+2z^{2} -45z^{4} + 228z^{5} - 230z^{6} + 2316z^{7} -
375z^{8} + 7020z^{9} + 296z^{10} + 5880z^{11} +
532z^{12}}{4z^{2}(1+9z^{2} +25z^{4} - 11z^{5} + 20z^{6} -
44z^{7})^{2}}$\\
\\
49Z \> \tiny $(1+4z+8z^{2} + 56z^{3} + 236z^{4} + 628z^{5} +
1346z^{6} + 992z^{7} - 1780z^{8} - 1880z^{9} + 1064z^{10} +
1100z^{11}$\\ \> \tiny $+925z^{12} - 1380z^{13} +
300z^{14})/(4z^{2}(1 + 2z - 14z^{2} - 35z^{3} + 14z^{4} + 49z^{5} -
21z^{6} - 11z^{7} + 6z^{8})^{2})$\\
\\
50Z \> \tiny $(1 + 12z + 72z^{2} + 360z^{3} + 1910z^{4} + 9488z^{5}
+ 38156z^{6} + 120216z^{7} + 300215z^{8} + 604300z^{9} +
990924z^{10}$\\ \> \tiny $+1326608z^{11} + 1441498z^{12} +
1256240z^{13} + 867180z^{14} + 476872z^{15} + 225869z^{16} +
111804z^{17} + 62000z^{18}$\\ \> \tiny $+30440z^{19} + 10436z^{20} +
2112z^{21} + 192z^{22})/(4z^{2}(-1 - 6z - 10z^{2} + 30z^{3} +
195z^{4} + 501z^{5} + 786z^{6} $\\ \>\tiny $+825z^{7} + 595z^{8} +
295z^{9} +
96z^{10} + 16z^{11})^{2})$\\
\end{tabbing}

\section{Transforming to any Constant}

These Q-values correspond to the modular function with a constant of
0.  To make our Q-values even more useful, we now give the details
on how to transform the Q-value to those of a modular function with
a different constant, say ``c."  After we've run our program, the
modular function is saved in \texttt{t[q]}.  Change \texttt{t[q]} to
\texttt{t[q]+c}.  Note - you must give a vlue for c.  Letting c be a
variable will make the calculations too complex for Mathematica to
handle.  Then run the remainder of the program, and the final
Q-value will be the Q-value determined by our modular function with
\texttt{c} as its constant.

\section{Bibliography}




\newdimen\invertedparindent
\invertedparindent=20pt

\newenvironment{invertedparagraphs}{%
\list{}%
{\setlength{\listparindent}{-\invertedparindent}%
\setlength{\leftmargin}{\invertedparindent}%
\setlength{\labelwidth}{0pt}%
\setlength{\labelsep}{0pt}%
\setlength{\itemindent}{0pt}}%
{\item{}\hspace{-\invertedparindent}}}{\endlist}
\begin{invertedparagraphs}

Abell, Martha L, and James P. Braselton.  \emph{Modern Differential
Equations}.  2nd ed.  New York:  Harcourt College Publishers, 2001.

Borcherds, Richard E.  "Monstrous moonshine and monstrous Lie
superalgebras."  Invent. Math. 109 (1992): 405 - 444.

Boyce, William E. and Richard C. DiPrima.  \emph{Elementary
Differential Equations and Boundary Value Problems}.  New York: John
Wiley and Sons, Inc., 1965.

Conway, J.H. and S.P. Norton.  "Monstrous Moonshine."  Bulletin of
the London Mathematical Society.  11 (1979): 308-339.

Hille, Einar.  \emph{Ordinary Differential Equations in the Complex
Domain}.  Mineola, NY: Dover Publications, Inc., 1976.

Lian, Bong H. and Shing-Tung Yau.  "Arithmetic Properties of Mirror
Map and Quantum Coupling."  Commun.Math.Phys.  176 (1991): 36 pp.
Online.  Internet: 13 Nov. 2002.

Lian, Bong H. and Shing-Tung Yau.  "Mirror Maps, Modular Relations
and Hypergeometric Series I."  Nuclear Phys. B Proc. Suppl. 46
(1996): 248 - 262.

Lopez, Robert J.  \emph{Advanced Engineering Mathematics}. Boston:
Addison Wesley, 2001.

Seaborn, James B.  \emph{Hypergeometric Functions and Their
Applications}.  New York: Springer-Verlag, 1991.

Zill, Dennis G.  \emph{A First Course in Differential Equations}.
Fifth Edition.  Boston: PWS-KENT Publishing Company, 1993.
\end{invertedparagraphs}


\end{document}